\documentclass[ 11pt] {article}
\usepackage{a4,  psfrag, epsfig, epsf, amscd, amsfonts, amsmath, amssymb, amstext, amsthm, enumerate}
\begin{document}

\newtheorem{satz}{Theorem}[section]
\newtheorem{defin}[satz]{Definition}
\newtheorem{hl}[satz]{Proposition}
\newtheorem{ko}[satz]{Corollary}
\newtheorem{que}[satz]{Question}
\newtheorem{bem}[satz]{Remark}

\title{A Remark to the Theorem of \\ Le Calvez and Yoccoz}
\author{Christian Pries}
\maketitle

\begin{abstract}
The theorem of Le Calvez and Yoccoz states that there are no minimal homeomorphisms on the finite punctered 2-dimensional sphere $S^{2}$. We show that this does not hold for other surfaces. Moreover, we discuss why the fast-conjugation-method fails in the most cases to construct such homeomorphism.
\end{abstract}

\setcounter{section}{0}
\section{Introduction}
It is well-known which closed 2-dimensional manifolds admit minimal homeomorphisms, but not for closed surfaces which are punctered. The case of the sphere was solved 1997:

\begin{satz}
(Le Calvez and Yoccoz) Given a finite set $F \subset S^{2}$. There is no minimal homeomorphism on $S^{2} - F$.
\end{satz}
Proof: (see [CY]) {\hfill $\Box$ \vspace{2mm}} \\ One can ask further if there are there counterexamples for other surfaces. We show that for two surfaces the answer is yes:
\begin{satz}
For any given non empty finite set $F$ of an orientable closed surface $X$ of genus equal $1$ or $2$, there is a minimal $C^{\infty}$ diffeomorphism on $X - F$.
\end{satz} 

\section{Construction Of Quasi-Minimal Systems \\ By Minimal Flows}

\begin{defin}
Let $(X,d)$ be a metric space and $G$ a topological group. A dynamical system $\Phi : G \times X \to X$ is called quasi-minimal if the union of its dense orbits is open. We call the non dense orbits the exceptional set.
\end{defin}

\begin{satz}
There exists a quasi-minimal flow on an orientable surface $X$ of genus $g = 2$ whose exceptional set is a point.
\end{satz} 

Proof: See section 14.4.b and section 14.6 in [KH] {\hfill $\Box$ \vspace{2mm}} \\ It is not clear if theorem [2.2] holds for surfaces of higher genus. The fixed point set for other quasi-minimal flows in [KH] corresponds the genus. For our examples we need an application of the following theorem:

\begin{satz}
Given a real topologically transitive flow $\Phi : \mathbb{R} \times X \to X$ on a separable metric space and suppose there is no "isolated streamline". For all values of $t \in \mathbb{R}$, except a set of first category, the homeomorphisms $f = \Phi_t: X \to X$ are topologically transitive.
\end{satz}

Proof: See [OU] {\hfill $\Box$ \vspace{2mm}}

\begin{ko}
Let $X$ be a manifold of dimension $n > 1$ and $\Phi : \mathbb{R} \times X \to X$ a $C^r$ quasi-minimal flow where $r > 0$. If its induced vectorfield is bounded, then for all values of $t \in \mathbb{R}$, except a set of first category, the $C^1$ diffeomorphisms $f = \Phi_t: X \to X$ are quasi-minimal and their exceptional sets coincide with the exceptional set of the flow $\Phi$.
\end{ko} Proof: Take a generic $t$ such that $f = \Phi_t: X \to X$ is transitive. Choose a $x_0$ such that $O_{f,+}(x_0)$ and $O_{f,-}(x_0)$ is dense. This is a well-known fact and follows from Baire's theorem. We show if any orbit $O_{\Phi,+}(x)$ is dense then $O_{f,+}(x)$ is dense too. Since the induced vectorfield is bounded, the orbit $O_{f,+}(x)$ must have an accumulation point on a segment of the orbit $O_{\Phi,+}(x_0)$. The accumulation point is dense, since $O_{f,+}(x_0)$ is dense and for any number $s$ we have $O_{f,+}(\Phi(s,x_0)) = \Phi(s,O_{f,+}(x_0))$, hence $O_{f,+}(x)$ is dense. The same argument works for the negative orbits. {\hfill $\Box$ \vspace{2mm}} 

Proof of theorem [1.2]: Due to thereoem [2.2] we know that if $X$ has genus 2 it admits at least a quasi-minimal flow with only one exceptional point. For genus 1 we take a minimal flow. Given a finite set $F \subset X$. We can assume that all elements of $F$ belongs to distinct orbits which are dense on both directions and the exceptional point belongs to $F$, otherwise we take a conjugation of the flow. We multiply the induced vector field $V$ of the flow with a function $ 0\le f$ that is exactly zero on $F$. The flow of $fV$ is quasi-minimal and the exceptional set is $X - F$. Now apply corollary [2.4] {\hfill $\Box$ \vspace{2mm}}

\section{The Fast-Conjugation-Method}

Fayad and Herman showed in [FH] that any compact manifold $M$ that admits a smooth free $S^1$ action must admit minimal diffeomorphisms. They used the fast-conjuagtion-method and the theorem of Baires to proof that in the closure of the set $\{ g\Phi_tg^{-1} | g\in $Diff($M$)$,  t\in S^1 \}$ with respect to the $C^r$ topology the minimal diffeomorphisms are generic. We can not extended this generic result to non compact manifolds or manifolds with a semi-free smooth $S^1$ action:
\begin{bem}
Given a continous $S^1$ action on a manifold $M$. For any subset $C \subset S^1$ and $G \subset$Diff($M$) the closure of $F_{C,G} = \{ g\Phi_tg^{-1} | g\in G,  t\in C \}$ contain a generic subset $D_{C,G}$ such that for each element in $D_{C,G}$ every orbit is postive recurrent,
hence any quasi-minimal homeomorphism in $D_{C,G}$ is minimal.
\end{bem}

Proof: Take two bases  $\{ U_i \}$ and $\{ V_i \}$ of the topology of $M$ such that $\overline U_i \subset V_i$ and set $R_i = \{ f \in$ Diff($M$) $ |$ $  \forall$ $ x \in \overline U_i$ $  \exists$ $ n_x > 0 : f^{n_x} \in V_i \}$. Notice that since the $S^1$ action is continous, we have $F_{C,G} \subset R_i$. Since the space Diff($M$) is a Baire space and the sets $R_i$ are open, the set $$D:= \overline F_{C,G} \cap \bigcap_{i} R_i $$ is generic in the closure of $\overline F_{C,G}$. We have that every point $x \in M$ is positive recurrent for any $f \in D$. Indeed, take a sequence of neighbourhoods $V_{j(i)}$ such that $\overline V_{j(i)} \to x$. If $x$ is not a periodic point, we can build at least a squence $n_i \to \infty$ such that $f^{n_i}(x) \to x$, so every point is positive recurrent. If $f \in D$ is quasi-minimal, then for every dense orbit we have $\overline O_{f,+}(x) = \overline O_{f}(x)$, since each orbit is positive recurrent, thus $f$ is a forward minimal homeomorphism on a dense open set. Due to theorem B of [G], there is no forward minimal homeomorphism on any non compact locally compact space, thus the open set of dense orbits is compact, so the homeomorphism is forward minimal. {\hfill $\Box$ \vspace{2mm}}

On the even-dimensional spheres $S^{2n}$ we can find semi-free smooth $S^1$ actions that are free except on two fixed points. Because of the last remark and the fact that $S^{2n}$ admits no minimal homeomorphism, we can not find a set $F_{C,G}$ with a generic subset of quasi-minimal diffeomorphism.

\vspace*{0,25 cm}

\end{document}